\documentclass[%
 amsmath,amssymb,
 prfluids,
]{revtex4-2}

\usepackage{newtxtext}
\usepackage{newtxmath}
\usepackage{natbib}
\usepackage{graphicx}
\usepackage{amsmath}
\usepackage{amssymb}
\usepackage{soul}
\usepackage{caption}
\usepackage{natbib}
\usepackage{soul}
\usepackage{color}
\usepackage{hyperref}
\hypersetup{
    colorlinks = true,
    urlcolor   = blue,
    citecolor  = black,
}
\newcommand{\RomanNumeralCaps}[1]

\definecolor{burgundy}{rgb}{0.5, 0.0, 0.13}

\definecolor{past}{rgb}{0.47,0.87,0.47} 

\def\bea{\begin{equation}}
\def\eea{\end{equation}}
\newcommand{\diffn}[3]{\frac{\mbox{d}^{#3} #1}{\mbox{d} #2^{#3}}}
\newcommand{\diff}[2]{\frac{\mbox{d} #1}{\mbox{d} #2}}
\newcommand{\pdiffn}[3]{\frac{\partial^{#3} #1}{\partial #2^{#3}}}
\newcommand{\pdiff}[2]{\frac{\partial #1}{\partial #2}}
\newcommand{\pdiffm}[3]{\frac{\partial^2 #1}{\partial #2 \partial #3}}

\def\diffoperator{\mathcal{D}}
\def\diffslowoperator{\mathcal{S}}
\def\shrod {Schr\"{o}dinger }

\begin{document}

\preprint{APS/123-QED}

\title{Parameter-free wave-envelope evolution equations including weak dissipation and forcing}

\author{J. S. Keeler}
\email{j.keeler@uea.ac.uk}
\affiliation{School of Engineering, Mathematics and Physics, University of East Anglia, Norwich, NR4 7TJ, UK}
\author{A. Alberello}
\email{a.alberello@uea.ac.uk}
\affiliation{School of Engineering, Mathematics and Physics, University of East Anglia, Norwich, NR4 7TJ, UK}
\author{B. Humphries}
\email{b.humphries@uea.ac.uk}
\affiliation{School of Engineering, Mathematics and Physics, University of East Anglia, Norwich, NR4 7TJ, UK}

\author{E. P{\u{a}}r{\u{a}}u}
\email{e.parau@uea.ac.uk}
\affiliation{School of Engineering, Mathematics and Physics, University of East Anglia, Norwich, NR4 7TJ, UK}

\begin{abstract}

  The higher-order nonlinear \shrod equation (Dysthe's equation in the context of water-waves) models the time evolution of the slowly modulated amplitude of a wave-packet in dispersive partial differential equations (PDE). These systems, of which water-waves are a canonical example, require the presence of a small-valued ordering parameter so that a multi-scale expansion can be performed. However, often the resulting system itself contains the small-ordering parameter. Thus, these models are difficult to interpret from a formal asymptotics perspective. This paper derives a parameter-free, higher-order evolution equation for a generic infinite-dimensional dispersive PDE with weak linear damping and/or forcing. Instead of focusing on the water-wave problem or another specific problem, our procedure avoids the complicated algebra by placing the PDE in an infinite-dimensional Hilbert space and Taylor expanding with Fre\'chet derivatives. An attractive feature of this procedure is that it can be used in many different physical settings, including water-waves, nonlinear optics and any dispersive system with weak dissipation or forcing. The paper concludes by discussing two specific examples.

\end{abstract}
\maketitle
\section{Introduction}

The nonlinear \shrod equation, first derived in \cite{zakharov1968stability}, and the higher-order version, or, in the context of water-waves, Dysthe's equation, first derived in \cite{dysthe1979}, are immensely successful equations. They describe the slowly-varying wave envelope of dispersive partial differential equations (PDE) in the field of water-waves \cite[for example]{lannes2013book} and optics \cite[for example]{boydoptics} and are cornerstone achievements of nonlinear science.

A key feature of the nonlinear \shrod equation (hereafter denoted NLS) and the higher-order \shrod equation (hereafter denoted HNLS) is that they are Hamiltonian systems (in the context of water-waves, for the HNLS this has only been shown very recently \cite{craig2021normal,guyenne2022vorticity}). In such systems, there is no dissipation mechanism. Therefore dynamical-attractors or repellers cannot exist and, crucially, energy cannot be gained or lost. In some physical situations, where energy-loss/gain is an important part of the physics, this is undesirable \cite{stuhlmeier2024modulational}. For example, in the evolution of ocean swell  \cite{segur2005stabilizing,henderson2013role} and wind forcing of ocean waves \cite{armaroli2018,eeltink2019,eeltink2020nonlineardyn} the dissipation results in a constant decay rate of the wave attenuation. Extending this, in ice-infested ocean-waves, experimental evidence suggests that the energy of a wave-packet is actually damped at a rate that depends on the wave frequency \cite{meylan2018jgr}. The classical NLS and HNLS, which in this context are derived from the inviscid, irrotational and incompressible free-surface Euler system \cite[for example]{lannes2013book}, cannot capture this energy attenuation. In order for the system to dissipate energy, an \textit{ad hoc} modification to the form of the NLS/HNLS is necessary as the Euler system is derived from the underlying Navier-Stokes system under the assumption of inviscid flow \cite{mei2003weakly,liao2023modified}. A rational argument of how the Navier-Stokes equations can lead to the free-surface Euler system with dissipation was proposed in \cite{dias2008pla} and a HNLS with dissipation was derived in \cite{carter2016frequency} but it is unclear how this analysis can be extended to more general dispersive systems with weak dissipation.

In early models of a modified-NLS \cite{mei2003weakly,segur2005stabilizing,henderson2013role} and modified-HNLS \cite{carter2016frequency}, an extra constant, linear, dissipation term was added. More recently, in light of striking evidence that the damping is frequency-dependent \cite{meylan2018jgr}, a non-constant, linear dissipation term was added to the NLS and HNLS to model sea-ice \cite{eeltink2020nonlineardyn,alberello2022dissipative,alberello2023dynamics}. In these models the equation is often simply stated without derivation and the coefficients often vary from study to study, even for identical physical situations (see for example, in the context of water-waves, \cite{OSBORNE2010573,peregrine1983water,zakharov1968stability,eeltink2019}. However, attempts to derive the equations from the fully nonlinear Euler system can be unwieldy, laborious, and prone to inadvertent algebraic mistakes; making derivations difficult to follow. In this paper, we derive a modified-HNLS system for a \textit{general} dispersive evolution PDE with weak dissipation. A key feature is that, in contrast to the literature, the small ordering parameter is absent from the final form of the resulting HNLS; we thus produce an asymptotically consistent form of this system. 

The literature on the NLS and HNLS is vast. However, a new aim of this paper is to develop and describe a friendly `user-guide' and systematic method to derive these evolution equations. An important feature of this paper is that this procedure can be applied to a wide range of dispersive evolution PDEs with weak dissipation or weak forcing. In this paper, we thus aim to describe the key steps in the formal asymptotic expansion as transparently as possible. We believe that this article will be a valuable tool for future derivations in other physical contexts and that an article of this type is long overdue in the nonlinear science community.

It is important to observe from the outset that in the literature, the nomenclature of higher-order NLS models and so-called `modified' NLS models have become convoluted (for example the modified NLS in \cite{liao2023modified} refers to the damping whereas in \cite{dysthe1979} it refers to the higher-order terms). To avoid any confusion, in this paper, we shall use the term modified-HNLS model, where `modified' refers to the dissipation term (similar for the NLS in \cite{liao2023modified}). 
\subsection*{Current approaches}

\begin{figure}
    \centering
    \includegraphics[width=\textwidth]{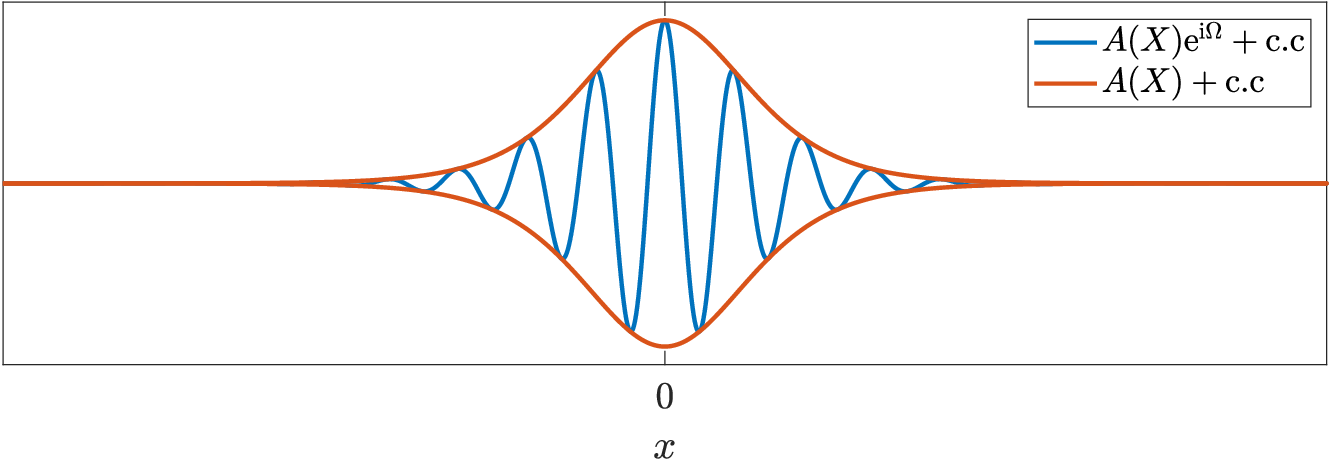}
    \caption{Sketch of a steadily translating wave packet, labelled $u={A}(X)\mbox{e}^{\mathrm{i}\Omega} + \mbox{c.c}$; $\Omega = {k}\boldsymbol{x}-\omega t$, with oscillations on a fast timescale (blue) and an envelope, ${A}(X)$, varying on a slow timescale (red).}
    \label{fig:wave_packet}
\end{figure}

Despite the success of the NLS (and associated systems) in describing the wave-amplitude for water-waves, from a formal asymptotic perspective, an undesirable aspect of the current models in the literature is the presence of a small, real parameter $\varepsilon\in\mathbb{R}$, where $|\varepsilon|\ll 1$. To illustrate this approach, let $u(x,t)$ represent the solution of a nonlinear PDE in $O(1)$ time and space variables, $x$ and $t$, written functionally as
\bea
\pdiff{u}{t} + \mbox{Nonlinear operator}(u;\varepsilon) = 0,
\label{pde12}
\eea
where the nonlinear operator contains spatial derivatives and is dependent on $\varepsilon$. Usually, to derive the NLS and HNLS systems (and modified versions) a solution to \eqref{pde12} is proposed of the form
\bea
u =A(\xi,\tau)\times \mbox{Periodic-function}(x,t) + \varepsilon\times\mbox{Higher-order terms}(x,t,\xi,\tau),
\label{usol}
\eea
where ${A}(\xi,\tau)$ is the complex-valued wave envelope depending on some slow space and time variables, $\xi$ and $\tau$, and the periodic-function is linear and periodic in both space and time, e.g. an exponential (see figure~\ref{fig:wave_packet}). For the modified HNLS the evolution of the wave-envelope, $A(\xi,\tau)$ is found by solving a nonlinear PDE of the form 
\begin{multline}
    \pdiff{A}{\tau} + \mbox{\shrod operator}(A) + \mbox{Dissipative operator}(A)  \\ + \varepsilon\times\left[\mbox{Higher-order \shrod operator}(A) + \pdiff{\:\mbox{Dissipative operator}(A)}{x}\right]= 0,
\label{dysthe_basic1}
\end{multline}
where the \shrod and higher-order \shrod operators contain linear spatial derivatives of $A$ and nonlinear terms involving $A$ and $\overline{A}$ (an over-bar indicates the complex conjugate) and the dissipation operator depends on the type of damping required \cite{eeltink2020nonlineardyn,alberello2022dissipative,alberello2023dynamics}. In these traditional approaches \eqref{dysthe_basic1} reduces to the `normal' nonlinear \shrod system, simply by setting $\varepsilon = 0$. 

We emphasise that models of the form in \eqref{dysthe_basic1} are extremely successful at modelling and predicting the evolution of the wave-packet \cite{segur2005stabilizing,carter2016frequency} yet, from a purely formal asymptotics perspective there is a troubling aspect of these equations in that the small-ordering parameter, $\varepsilon$, is present in the system as a coefficient of the higher-order \shrod operator. Furthermore, from a practical perspective, for numerical simulations, $\varepsilon$ is a control parameter, rather than an ordering parameter, which can result in unnecessary `stiffness' as $\varepsilon\ll 1$ in the numerical implementation. Therefore, the aim of this article is to describe definite procedures for deriving \textit{parameter-free}, higher-order evolution equations for a general dispersive system with weak dissipation.

\subsection*{A parameter-free approach}

We now sketch our general strategy to achieve a parameter-free modified HNLS. One way to achieve this is by introducing a succession of small time-scales proportional to $\varepsilon^n$ \cite{degasperis1997multiple,calogero2000nonlinear}. With the view of a numerical implementation, where the discretisation of different time-scales may be problematic, we will expand the amplitude function instead, resulting in a coupled system of evolution equations. These are more amenable to numerical study as there is only one scale for the independent variable in the slow-time scale, $\tau$. Our strategy involves expanding the amplitude function to higher-orders so that instead of \eqref{usol} we expand the solution using 
\bea
u =\left[A(\xi,\tau) + \varepsilon B(\xi,\tau)\cdots\right]\times \mbox{Periodic-function}(x,t) + \varepsilon\times\mbox{Higher-order terms}(x,t,\xi,\tau),
\eea
where ${B}(\xi,\tau)$ is the slowly-varying first-order envelope function. We will show that an asymptotically consistent system for ${A}$ and ${B}$, which is \textit{parameter-free} is of the form
\begin{align}
&\pdiff{A}{\tau} + \mbox{Linear derivative} + \mbox{Nonlinear term} + \mbox{Dissipation term} = 0,\\
&\pdiff{B}{\tau} + \mbox{Linear derivative} + \mbox{Nonlinear terms} + \mbox{Dissipation term}= 0.
\end{align}

This approach relies on the interpretation of an evolution PDE as an infinite dimensional dynamical system where, instead of state variables belonging to $\mathbb{R}^d$, with $d$ finite, the state variables, $u$, belong to an infinite-dimensional Hilbert space and is, in fact, inspired by work where a similar formal asymptotic expansion was used in a problem involving the propagation of air-bubbles in a viscous fluid within a Hele-Shaw channel \cite{keeler2019invariant}. 

We proceed as follows. In \S~\ref{sec:prem}, we will devote the preliminary part of this article to describing this framework. Then in \S~\ref{sec:formal} we will proceed with a formal asymptotic analysis. First, we will show how a naive asymptotic expansion fails and then provide an appropriate remedy. In our analysis, our main philosophy is to be constructive and transparent to the reader, so that this article be a `user-friendly' guide. In \S~\ref{sec:kdv} and \ref{sec:water}, we focus on the derivation from a toy system and the fully nonlinear water-wave problem, respectively. For the latter, the rationale for this is that although the evolution amplitude equations for water-waves are well known in the case of simple geometries (infinite depth or flat finite depth) and for `standard' physics (i.e. gravity-capillary waves) \cite{lannes2013book}, a researcher may wish to determine the evolution equations for their particular problem which may have different geometries (for example a submerged body in the fluid, or a variable bottom topography) and include different physics (for example hydroelastic waves, external pressure distributions). Therefore the secondary aim of this article is to provide a useful user-friendly `look-up' tool for a researcher so they can make the necessary adjustments for their own physical model. 

\section{Preliminaries}\label{sec:prem}
We will study the generic PDE for the evolution of $u$:
\bea
\pdiff{u}{t} + \mathcal{F}(u;\varepsilon) + \varepsilon^2\mathcal{V}(u) = 0,\qquad x\in\mathbb{R}^d,t\in\mathbb{R}^+,
\label{pde}
\eea
where $u(x,t)\in\mathcal{U}$, a Hilbert space, $\mathcal{F}:\mathcal{U}\mapsto \mathcal{U}$ is a nonlinear operator that depends on a parameter $\varepsilon \in \mathbb{R}$ and $\mathcal{V}:\mathcal{U}\mapsto \mathcal{U}$ is a linear dissipation operator. We note that $\mathcal{F}(u;\varepsilon)$ only contain spatial derivatives and in the rest of the paper, for ease of exposition, will omit the explicit $\varepsilon$ dependence from $\mathcal{F}$. We also assume that associated with \eqref{pde} are the correct number of boundary conditions on $u$ and spatial derivatives of $u$ to make the system well-posed. We assume that $d=1$ (this will be relaxed in \S~\ref{sec:water}) and that when $\varepsilon = 0$, the linearised equation admits a real-valued dispersion relation. We start our analysis with several definitions involving Hilbert spaces and Fr\'echet derivatives and introduce a notation that simplifies the subsequent analysis. First, we assume that $u = u_s(x)$ is a steady state of \eqref{pde}, independent of time, i.e.
\bea
\mathcal{F}(u_s) = 0.
\eea
We emphasise that $u_s$ does not have to be identically zero; the analysis will hold for a generic steady base state $u_s$. Anticipating periodic solutions of \eqref{pde} with temporal period $\tilde{T}$ and spatial period $\tilde{X}$, we  define an inner product
\bea
\langle u,v \rangle = \int_{0}^{\tilde{T}}\int_{0}^{\tilde{X}}u\overline{v}\,\mbox{d}x\,\mbox{d}t, \qquad u,v\in\mathcal{U}.
\eea

\subsection{Taylor expansions in Hilbert spaces}\label{sec:taylor}
An important feature of the subsequent analysis is the ability to Taylor expand the fully nonlinear operator, $\mathcal{F}$ within a Hilbert space  (see, for example \S~5.6 of \cite{cartan1971differential}). We assume that $\mathcal{F}(u)$ is $n$-differentiable, in the Fre\'chet sense (c.f. \eqref{frechet}), and admits a Taylor expansion about the steady base state, $u = u_s + v$, i.e.
\bea
\mathcal{F}(u) = \underbrace{\mathcal{F}(u_s)}_{=0} + D^1[\mathcal{F}(u_s)](v) + \frac{1}{2!}D^2[\mathcal{F}(u_s)](v,v) + \frac{1}{3!}D^3[\mathcal{F}(u_s)](v,v,v) + \cdots;
\label{taylor}
\eea
where the first-order Fre\'chet derivative is defined as
\bea
D^1[\mathcal{F}(u_s)](v_1) \equiv \lim_{h\to 0}\frac{\mathcal{F}(u_s + hv_1) - \mathcal{F}(u_s)}{h}.
\label{frechet}
\eea
Higher-order derivatives can be defined inductively;
\bea
D^n[\mathcal{F}(u_s)](v_1,\cdots,v_n)\equiv D[D^{n-1}(v_1,\cdots,v_{n-1})](v_n),\qquad v_i\in\mathcal{U}.
\eea
We emphasise that the $n$\textsuperscript{th} Fre\'chet derivative is a $n$-linear functional. For ease of exposition, we will introduce specific symbols for the first four Fre\'chet derivatives about $u_s$ as
\bea
\begin{split}
&\mathcal{J}(u_1) \equiv D^1[\mathcal{F}(u_s)](u_1),\quad \mathcal{H}(v_1,v_2)\equiv \frac{1}{2!}D^2[\mathcal{F}(u_s)](v_1,v_2),\\
&\mathcal{T}(v_1,v_2,v_3)\equiv \frac{1}{3!}D^3[\mathcal{F}(u_s)](v_1,v_2,v_3)  ,\quad\mathcal{Q}(v_1,v_2,v_3,v_4)\equiv  \frac{1}{4!}D^4[\mathcal{F}(u_s)](v_1,v_2,v_3,v_4).
\end{split}
\label{not}
\eea
We emphasise that the Taylor expansion in \eqref{taylor} is for a general nonlinear operator, $\mathcal{F}(u)$.  

\subsection{The linear operator}
We now discuss the linear operator, $\mathcal{J}(u)$ in \eqref{not}, that will perform an important role in the subsequent analysis. In what follows we adopt the notation
\bea
\Omega_i = k_ix - \omega t,\qquad \Omega_{i,j} = (k_i+k_j)x - \omega t,\qquad \Omega_{i,j,l} = (k_i+k_j+k_l)x - \omega t,
\label{omega_def}
\eea
where $k_i$ is a spatial wave number, $\omega$ is a temporal wave number or angular frequency and $i$ is a labelling index. We can map our operator onto Fourier space as a pseudo-operator, i.e.
\bea
\mathcal{J}(A\mbox{e}^{\mathrm{i}\Omega_1}) \mapsto A\mbox{e}^{\mathrm{i}\Omega_1}J(k_1),
\label{linear1}
\eea
where $J:\mathbb{R}\mapsto\mathbb{R}$ is a function of a real variable. 

\subsection{The bilinear and trilinear operators}
In a similar way, we can map the bilinear and trilinear operators, $\mathcal{H}$, $\mathcal{T}$, in Fourier space so that
\begin{align}
&\mathcal{H}(A_1\mbox{e}^{\mathrm{i}\Omega_1},A_2\mbox{e}^{\mathrm{i}\Omega_2}) \mapsto A_1A_2\mbox{e}^{\mathrm{i}\Omega_{1,2}}H(k_1,k_2),\label{bilinear}\\
&\mathcal{T}(A_1\mbox{e}^{\mathrm{i}\Omega_1},A_2\mbox{e}^{\mathrm{i}\Omega_2},A_3\mbox{e}^{\mathrm{i}\Omega_3}) \mapsto A_1A_2A_3\mbox{e}^{\mathrm{i}\Omega_{1,2,3}} T(k_1,k_2,k_3).\label{trilinear}
\end{align}
The bilinearity of $\mathcal{H}$ means the following identities hold
\bea
\mathcal{H}(a+b,c+d) \equiv\mathcal{H}(a,c) + \mathcal{H}(a,d) + \mathcal{H}(b,c) + \mathcal{H}(b,d),\qquad \mathcal{H}(\lambda a,\mu b) \equiv \lambda\mu\mathcal{H}(a,b),
\label{identities}
\eea
$\forall a,b,c,d\in\mathcal{U}$ and $\forall\lambda,\mu\in\mathbb{C}$. Similar distributive properties hold for the trilinear operator, $\mathcal{T}$. 

\section{Formal asymptotic expansion}\label{sec:formal}

Now we have defined our problem and the associated operators we can attempt to construct an asymptotically valid expansion for the solution. First, we perform a naive asymptotic expansion to solve \eqref{pde} and show how it fails, to understand how we may find a remedy. We now assume that $\varepsilon\ll 1$, and formally expand $u$ in the form
\bea
u(x,t) = u_s(x) + \varepsilon u_1(x,t) + \varepsilon^2 u_2(x,t) + \varepsilon^3 u_3(x,t) + \varepsilon^4 u_4(x,t) + \cdots,\qquad |\varepsilon|\ll 1.
\label{exp}
\eea
Substituting \eqref{exp} into \eqref{pde} and using \eqref{taylor} gives
\begin{multline}
    \varepsilon\pdiff{u_1}{t} + \varepsilon^2\pdiff{u_2}{t} + \varepsilon^3\pdiff{u_3}{t} + \cdots \underbrace{\left[\varepsilon \mathcal{J}(u^*) + \varepsilon^2\mathcal{H}( u^*,u^*)
+ \varepsilon^3\mathcal{T}( u^*,u^*,u^*) +\cdots \right]}_{\mathrm{Taylor\: expansion}} = 0,
\label{pde1}
\end{multline}
where we explicitly highlight the higher-order terms resulting from the Taylor expansion, see \eqref{usol}, and $u^* = u_1 + \varepsilon u_2 + \varepsilon^2u_3 + \cdots$.

\subsection{Leading-order: The dispersion relation}

In \eqref{pde1}, at $O(\varepsilon)$, we have
\bea
\mathcal{L}(u_1) = 0,\qquad\qquad \left[\mathcal{L}(\star)\equiv \pdiff{(\star)}{t} + \mathcal{J}(\star),\qquad \star \in \mathcal{U}\right].
\label{ordereps}
\eea
Equation \eqref{ordereps} has a solution
\bea
u_1(x,t) = A\mbox{e}^{\mathrm{i}\Omega} + \overline{A}\mbox{e}^{-\mathrm{i}\Omega},\qquad \Omega = kx - \omega t,
\label{u1}
\eea
for an, as yet, unknown constant $A\in\mathbb{C}$. Using our notation introduced in \eqref{linear1}, $k$ and $\omega$ are related through the linear dispersion relation
\bea
\omega + {J}(k)=0.
\label{disp}
\eea

\subsection{Disordered expansion}

Continuing our naive expansion of \eqref{pde1}, at $O(\varepsilon^2)$ we obtain
\bea
\mathcal{L}(u_2) = -\mathcal{H}(u_1,u_1).
\eea
The operator $\mathcal{H}$ is bilinear and examining the form of \eqref{u1} means that $\mathcal{H}$ will contain two types of terms proportional to i) $\mbox{e}^{2\mathrm{i}\Omega}$ and ii) $\mbox{e}^{0\mathrm{i}\Omega}$. Therefore the solution at this order will be of the form
\bea
u_2(x,t) = \varphi_0A^2\mbox{e}^{2\mathrm{i}\Omega} + B\mbox{e}^{\mathrm{i}\Omega} + \varphi_1|A|^2 + \mbox{c.c},
\eea
where $\varphi_{0,1}$ are known solutions to linear problems, $B\in\mathbb{C}$ is arbitrary and c.c. stands for complex conjugate. At this stage, $u_1$ and $u_2$ are bounded and periodic. At the next order $O(\varepsilon^3)$ we find
\bea
\mathcal{L}(u_3) = -\mathcal{H}(u_1,u_2) -\mathcal{H}(u_1,u_2) - \mathcal{T}(u_1,u_1,u_1),
\eea
Now there are terms on the right-hand side that are proportional to $A|A|^2\mbox{e}^{\mathrm{i}\Omega}$ and hence $u_3$ will experience secular growth unless we choose $A=0$, which reduces the expansion to null and thus illustrates the failure of the naive expansion.

\subsection{Multiple-scales and two-timing}\label{sec:twotimeexpansion}

To avoid a trivial selection of our amplitude function, we introduce slow variables $X$,$T$ that are related to the $O(1)$ quantities by
\bea
X = \varepsilon x,\qquad T = \varepsilon t.
\label{smalltime}
\eea
With $u(x,t)\mapsto u(x,t,X,T)$, the differential operators get mapped to
\bea
\pdiff{}{x}\mapsto \pdiff{}{x} + \varepsilon\pdiff{}{X}.
\eea
We have to be careful with how this ansatz modifies the multi-linear operators introduced in the previous section. Due to this ansatz, the wave number in Fourier space is mapped to
\bea
k\mapsto k + \varepsilon K,
\label{kmap}
\eea
where $K$ is a slow wave number associated with the slow spatial scale, $X$. Using \eqref{kmap}, expanding $\mathcal{J}$ as a Taylor series yields
\begin{align}
\mathcal{J}(A(X,T)\mbox{e}^{\mathrm{i}\Omega}) &\mapsto \left[J(k) + \varepsilon \diff{J}{k}K + \frac{1}{2}\varepsilon^2\diffn{J}{k}{2}K^2 + \cdots\right]A(X,T)\mbox{e}^{\mathrm{i}\Omega},\\
&\mapsto \left[J(k) + \varepsilon \diff{J}{k}\pdiff{A}{X} + \frac{1}{2}\varepsilon^2\diffn{J}{k}{2}\pdiffn{A}{X}{2} + \cdots\right]\mbox{e}^{\mathrm{i}\Omega}.\label{asd2}
\end{align} 
For the bilinear operator, we have to be more careful. The wave numbers are mapped according to $k_1\mapsto k_1 + \varepsilon K_1,k_2\mapsto k_2+\varepsilon K_2$ so that
\begin{align}
    \mathcal{H}(A_1\mbox{e}^{\mathrm{i}\Omega_1},A_2\mbox{e}^{\mathrm{i}\Omega_2}) &\mapsto \left[H(k_1,k_2) + \varepsilon\boldsymbol{k}^T\nabla_{k}H+\frac{1}{2}\varepsilon^2\boldsymbol{K}^T\boldsymbol{H}\boldsymbol{K}+\cdots \right]A_1A_2\mbox{e}^{\mathrm{i}\Omega_{1,2}},\label{asd}\\
    &\mapsto\left[H(k_1,k_2) + \varepsilon\left(\pdiff{H}{k_1}\pdiff{A_1}{X}A_2 + \pdiff{H}{k_2}A_1\pdiff{A_2}{X}\right) + \right.\label{asd1}\\
&\left.\frac{1}{2}\varepsilon^2\left(\pdiffn{H}{k_1}{2}\pdiffn{A_1}{X}{2}A_2 + 2\pdiffm{H}{k_1}{k_2}\pdiff{A_1}{X}\pdiff{A_2}{X}+\pdiffn{H}{k_2}{2}A_1\pdiffn{A_2}{X}{2}\right)+\cdots\right]\mbox{e}^{\mathrm{i}\Omega_{1,2}},
\nonumber\end{align}
where the slow wave number vector is defined as $\boldsymbol{K} = (K_1,K_2)^T$, wave number gradient $\nabla_k = (\partial/\partial k_1,\partial/\partial k_2)$  and $\boldsymbol{H}$ is the Hessian matrix of $H$ with respect to $k_1,k_2$. When going from \eqref{asd} to \eqref{asd1} we consider that in Fourier space $A_1 \mapsto \hat{A}_1(K_1,T)$ and $A_2 \mapsto \hat{A}_2(K_2,T)$.

Finally, we do the same to the trilinear operator, this time the slow wave number vector is $\boldsymbol{K} = (K_1,K_2,K_3)^T$ and gradient operator, $\nabla_k = (\partial /\partial k_1,\partial/\partial k_2,\partial/\partial k_3)^T$, where the ambiguity with the expressions above is noted, yet accepted as a reasonable abuse of notation. The trilinear operator is mapped to
\begin{align}
    \mathcal{T}(A_1\mbox{e}^{\mathrm{i}\Omega_1},A_2\mbox{e}^{\mathrm{i}\Omega_2},A_3\mbox{e}^{\mathrm{i}\Omega_3}) &\mapsto \left[T(k_1,k_2,k_3) + \varepsilon\boldsymbol{K}^T\nabla_k T+\cdots\right]A_1A_2A_3\mbox{e}^{\mathrm{i}\Omega_{1,2,3}}.\label{asd3}
\end{align}
It will be useful later to expand the gradient term above explicitly:
\bea
  \left[\boldsymbol{K}^T\nabla_k T\right]A_1A_2A_3 = \pdiff{T}{k_1}\pdiff{A_{1}}{X}A_2A_3+\pdiff{T}{k_1}A_1\pdiff{A_2}{X}A_3+\pdiff{T}{k_1}A_1A_2\pdiff{A_3}{X}.
\eea
For ease of exposition, we make the following definitions: 
\begin{equation}
\begin{split}
\diffslowoperator^n_{\mathcal{J}}(u_1) \equiv \frac{1}{n!}\diffn{J}{k_1}{n}\pdiff{A}{X}\mbox{e}^{\mathrm{i}\Omega},\qquad \diffslowoperator_{\mathcal{H}}(u_1,u_2) \equiv \left[\boldsymbol{K}^T\nabla_kH\right]A_1A_2\mbox{e}^{\mathrm{i}\Omega_{1,2}}\\
\diffslowoperator^2_{\mathcal{H}}(u_1,u_2) \equiv \frac{1}{2!}\left[\boldsymbol{K}^T\boldsymbol{H}\boldsymbol{K}\right]A_1A_2\mbox{e}^{\mathrm{i}\Omega_{1,2}}, \diffslowoperator_{\mathcal{T}}(u_1,u_2,u_3) \equiv \frac{1}{3!}\left[\boldsymbol{K}^T\nabla_kT\right]A_1A_2A_3\mbox{e}^{\mathrm{i}\Omega_{1,2,3}}.
\end{split}
\label{importantnotation}
\end{equation}
We have used the symbol $\mathcal{S}$ to denote an operator defined on the `slow' scale. For a technical discussion of the multi-scale ansatz see \cite{craig1992nonlinear}.

\subsection{The two-timed system}
Using the notation introduced in \S~\ref{sec:twotimeexpansion}, we can now write our system in a way that will let us to pick out the terms at each order in a straightforward manner. 

Summarising the process so far, we have performed two sets of Taylor expansions. The first involved a Hilbert space Taylor-expansion of the fully nonlinear operator $\mathcal{F}$ about the base state, $u_s$. This procedure resulted in a succession of multi-linear operators, $\mathcal{J},\mathcal{H}, \dots$, at each order in $\varepsilon$. The second set of Taylor expansions involved expanding the individual multi-linear operators due to the multi-scale ansatz in \eqref{smalltime}.

Bringing these expansions all together, using the notation in \eqref{importantnotation}, we can write \eqref{pde1} to an arbitrary order in $\varepsilon$. The new expansion is of the form
\bea
u = u_s + \varepsilon A(X,T)\mbox{e}^{\mathrm{i}\Omega} + \varepsilon^2 u_2(X,T,x,t) + \varepsilon^3 u_3(X,T,x,t) + \varepsilon^4 u_4(X,T,x,t) + \cdots.
\label{usol_best}
\eea
Substituting \eqref{usol_best} into \eqref{pde1}, and utilising the Taylor expansions in \eqref{importantnotation}, yields
\begin{multline}
  \underbrace{\varepsilon\left(\pdiff{u_1}{t}+\varepsilon\pdiff{u_1}{T}\right) + \varepsilon^2\left(\pdiff{u_2}{t}  +\varepsilon\pdiff{u_2}{T}\right)+ \varepsilon^3\left(\pdiff{u_3}{t}+\varepsilon\pdiff{u_3}{T}\right)}_{\mathrm{Time\, derivatives}} \\ +\underbrace{\varepsilon\mathcal{J}(u^*) + \varepsilon^2\diffslowoperator_{\mathcal{J}}(u^*) + \varepsilon^3\diffslowoperator_{\mathcal{J}}^2(u^*) + \varepsilon^4\diffslowoperator_{\mathcal{J}}^3(u^*)}_{\mathrm{Linear\:dispersion},\: \mathcal{J}}\\
   +\underbrace{\varepsilon^2\mathcal{H}(u^*,u^*) + \varepsilon^3\diffslowoperator_{\mathcal{H}}(u^*,u^*) + \varepsilon^4\diffslowoperator_{\mathcal{H}}^2(u^*,u^*)}_{\mathrm{Bilinear\:terms},\:\mathcal{H}}\\+
    \underbrace{\varepsilon^3\mathcal{T}(u^*,u^*,u^*) + \varepsilon^4\diffslowoperator_{\mathcal{T}}(u^*,u^*,u^*) }_{\mathrm{Trilinear\:terms},\:\mathcal{T}} +  \underbrace{\varepsilon^4\mathcal{Q}(u^*,u^*,u^*,u^*)}_{\mathrm{Quadrilinear\: term}}+ \underbrace{\varepsilon^2\mathcal{V}(u^*)}_{\mathrm{Dissipation}} + O(\varepsilon^5)= 0.
\label{pde5}
\end{multline}
Each set of grouped terms corresponds to the Taylor expansion of each of the (multilinear) operators in \eqref{pde1} due to the two-timing ansatz \eqref{smalltime}. 

The system in \eqref{pde5} may seem unwieldy. Still, conceptually it is simple as it arises from Taylor expansions in Hilbert spaces, which have a strong correlation to `standard' Taylor expansions of functions of real variables, with which we assume the reader is familiar. In practice all that is required are the i) elementary linear dispersion relation and ii) the forms of the bilinear and trilinear operators. For a given system the linear dispersion relation is often a trivial exercise, but the form of the bilinear and trilinear operators often requires some work. However, once established, the subsequent analysis becomes straightforward as it simply involves repeated applications of these operators with different arguments. 

\subsection{Two-time first order}
We now proceed with the formal asymptotic expansion. As before, at $O(\varepsilon)$, we obtain 
\bea
\mathcal{L}(u_1) = 0,
\eea
but this the time solution for $u_1$ will contain an unknown slowly varying (in space and time) modulation function;
\bea
u_1(x,t,X,T) = \underbrace{A(X,T)\mbox{e}^{\mathrm{i}\Omega} + \overline{A}(X,T)\mbox{e}^{-\mathrm{i}\Omega}}_{\mathrm{Oscillating}} ,\qquad \Omega = kx - \omega t\qquad (\mbox{\textbf{Leading-order solution}}), 
\eea
with the dispersion relation stated in \eqref{disp}.

\subsection{Two-time second order: The transport equation}

We continue our expansion to $O(\varepsilon^2)$ where
\bea
\mathcal{L}(u_2) + \pdiff{u_1}{T} + \diffslowoperator_{\mathcal{J}}(u_1) + \mathcal{H}(u_1,u_1) = 0,
\eea
or
\bea
\mathcal{L}(u_2) = -\underbrace{\left[\pdiff{A}{T} + c_g\pdiff{A}{X}\right]}_{\mathrm{Resonant}}\mbox{e}^{\mathrm{i}\Omega} - \underbrace{\mathcal{H}(u_1,u_1)}_{\mathrm{Non-resonant}} + \mbox{c.c},
\label{second}
\eea
where $c_g=\mbox{d}J/\mbox{d}k$ is the group velocity of the travelling wave. To eliminate the resonant terms on the right-hand side and thus avoid secular growth of the harmonic modes for $u_2$, we invoke the Fredholm alternative, see, \cite{griffel2002applied}. This states that either i) the only solution of $\mathcal{L}(u) = 0$ is $u=0$ and then $\mathcal{L}(u) = v$ has a unique solution \textit{or} ii) there exists a non-zero solution of $\mathcal{L}(u) =0$ in which case $\mathcal{L}(u)= v$ has a solution only if the inner product of $v$ with the solutions of $\mathcal{L}^\dagger(u^\dagger) = 0$ vanishes, where $\mathcal{L}^\dagger,u^{\dagger}$ are the adjoint problem and eigenmodes, respectively, of 
\bea
\langle \mathcal{L}(u),v\rangle = \langle u,\mathcal{L}^{\dagger}(v)\rangle,\qquad \forall u,v\in\mathcal{U}.
\eea
In what follows we assume our problem is self-adjoint but a non self-adjoint problem can be handled similarly \cite{keeler2019invariant}. For self-adjoint operators the adjoint eigenmodes are proportional to $\mbox{e}^{\mathrm{i}\Omega}$ and hence in our problem, \eqref{second}, this can be written down as
\bea
\langle\left(\pdiff{A}{T} + c_g\pdiff{A}{X}\right)\mbox{e}^{\mathrm{i\Omega}}+\mathcal{H}(u_1,u_1),u^{\dagger}\rangle = 0.
\label{transport1}
\eea
The non-resonant terms in $\mathcal{H}(u_1,u_1)$ will automatically be orthogonal to $u^\dagger$ due to periodicity and so the solvability condition in \eqref{transport1} becomes
\bea
\pdiff{A}{T} + c_g\pdiff{A}{X} = 0.
\label{transport}
\eea
We can write an alternative form of \eqref{second} using the Fourier representation,
\begin{multline}
  \mathcal{L}(u_2)= -\underbrace{A^2{H}(k,k)\mbox{e}^{2\mathrm{i}\Omega}  -\overline{A}^2{H}(-k,-k)\mbox{e}^{-2\mathrm{i}\Omega}}_{\mathrm{Oscillating}}-\underbrace{|A|^2H(k,-k) -|A|^2H(-k,k)}_{\mathrm{Non-oscillating}} \\
  + \mbox{Resonant terms}.
\end{multline}
This form allows us to write down the solution, $u_2$ as
\bea
u_2 =  \varphi_0\underbrace{A^2\mbox{e}^{2\mathrm{i}\Omega} +  \overline{\varphi}_0\overline{A}^2\mbox{e}^{-2\mathrm{i}\Omega}+B\mbox{e}^{\mathrm{i}\Omega} + \overline{B}\mbox{e}^{-\mathrm{i}\Omega}}_{\mathrm{Oscillating}} + \underbrace{{\varphi_{1}\mathbb{G}(|A|^2)}}_{\mathrm{Non-oscillating}}\qquad (\mbox{\textbf{First-order solution}}), 
\label{u2}
\eea
where $B(X,T)$ is an arbitrary modulation that can be considered as an order $\varepsilon$ correction to $A$ and $\varphi_i$ are functions of $k$ that can, in principle, be calculated. We believe that rather than providing explicit formulae for these coefficients, it is more instructive to demonstrate their determination through an explicit example, as we do in \S~\ref{sec:kdv} and \S~\ref{water_wave_dyn}. Finally, the operator $\mathbb{G}$ is defined as 
\bea
\mathbb{G}(u)=\int_{s\in\mathbb{R}^d} G(x,s)\,u\,\mbox{d}s,
\eea
where $G$ is the Green's function of $\mathcal{L}$, which we assume exists.

\subsection{Two-time third order: The modified-NLS system}

Continuing our analysis, at $O(\varepsilon^3)$ we have
\begin{multline}
\mathcal{L}(u_3) = -\left(\pdiff{u_2}{T}+c_g\pdiff{u_2}{X} + \frac{\mathrm{i}}{2}\diffn{\omega}{k}{2}\pdiffn{u_1}{X}{2}\right)  -\mathcal{H}(u_1,u_2) - \mathcal{H}(u_2,u_1)  \\ 
-\diffslowoperator_{\mathcal{H}}(u_1,u_1) - \mathcal{T}(u_1,u_1,u_1) - \mathcal{V}(u_1).
\label{pde2a}
\end{multline}
or in Fourier representation
\begin{multline}
    \mathcal{L}(u_3) = -\underbrace{\left[\pdiff{u_2}{T}+c_g\pdiff{u_2}{X}\right]}_{\mbox{Transport}} -\underbrace{\left[ \frac{\mathrm{i}}{2}\diffn{\omega}{k}{2}\pdiffn{A}{X}{2} + \lambda_1 A|A|^2 + V(k)A + {\lambda_2\mathbb{G}(|A|^2)}\right]\mbox{e}^{\mathrm{i}\Omega}}_{\mbox{Resonant}} \\
    -\underbrace{\left[\lambda_3A^3\mbox{e}^{3\mathrm{i}\Omega} + \left(\lambda_4A\pdiff{A}{X} + \lambda_5AB\right)\mbox{e}^{2\mathrm{i}\Omega} + \lambda_6\pdiff{|A|^2}{X} +\lambda_7(\overline{A}B+A\overline{B})\right]}_{\mathrm{Non-resonant}} + \mbox{c.c}
\label{pde3}
\end{multline}
where $\lambda_i$ are functions of $k$ that again, in principle, can be determined, and $V(k)$ is the Fourier symbol of the dissipation term. We note that the second term in the bracket containing the resonant terms arises due to the nonlinear interaction between $u_1$ and $u_2$ in those terms involving $\mathcal{H}$ and $\mathcal{T}$ in \eqref{pde2a}. A further important point is that a non-local term appears in the resonant part of the right-hand side of \eqref{pde3}. We note that non-local terms have appeared before in NLS-type models, for example, water-waves \cite{peregrine1983water}, a `shallow-deep' approximation for stratified fluids \cite{pelinovsky1996non} and in so-called generalised quasi-nonlinear models of shear flow \cite{marston2016generalized}. 

The transport term will vanish if we move to the frame of reference moving with the group velocity of the wave which the form of \eqref{transport} suggests. Therefore, we introduce a travelling wave coordinate, $\xi$, moving with speed $c_g$. In addition, the resonant term in \eqref{pde3} has zero time-derivative (in $T$), implying that $A$ is constant in time. To ensure a slow modulation, in time, of the envelope, $A$, we introduce a further slow time-scale, $\tau$, so that
\bea
\xi = X - c_gT,\qquad \tau = \varepsilon^2 t,\qquad A(X,T)\mapsto A(\xi,\tau),\qquad B(X,T)\mapsto B(\xi,\tau).
\label{travellingwave}
\eea
Therefore, \eqref{pde3} becomes
\bea
\mathcal{L}(u_3) = -\mbox{e}^{\mathrm{i}\Omega}\left(\pdiff{A}{\tau} + \frac{\mathrm{i}}{2}\diffn{\omega}{k}{2}\pdiffn{A}{\xi}{2} + \lambda_1 A|A|^2 + V(k)A + {\lambda_5A\mathbb{G}(|A|^2)}\right) + \mbox{Non-resonant terms}.
\label{pde4}
\eea
Now, to avoid secular growth, we repeat the imposition of the Fredholm alternative. Again, assuming $\mathcal{L}$ is self-adjoint, the solvability condition is the \textbf{modified-NLS equation:} 
\bea
\boxed{\pdiff{A}{\tau} + \frac{\mathrm{i}}{2}\diffn{\omega}{k}{2}\pdiffn{A}{\xi}{2} + \lambda_1 A|A|^2 + V(k)A + {\lambda_5A\mathbb{G}(|A|^2)} = 0.}
\label{A_eqn}
\eea
We emphasise that this equation is valid for any system of the form \eqref{pde} with small linear dissipation. We note that the non-local operator, $\mathbb{G}$, in the context of water-waves is absent as we will show in \S~\ref{sec:water}.

As in the previous order, using Fourier transforms, the full particular integral for $u_3$ can be written as
\bea
\begin{split}
&u_3 = \underbrace{A^3\psi_0\mbox{e}^{3\mathrm{i}\Omega} +\overline{A}^3\overline{\psi}_0\mbox{e}^{-3\mathrm{i}\Omega} + \left(\psi_1A\pdiff{A}{\xi}+\psi_2AB\right)\mbox{e}^{2\mathrm{i}\Omega} + \overline{\psi}_3\overline{A}\pdiff{\overline{A}}{\xi}\mbox{e}^{-2\mathrm{i}\Omega}+C\mbox{e}^{\mathrm{i}\Omega} + \overline{C}\mbox{e}^{-\mathrm{i}\Omega}}_{\mathrm{Oscillating}} \\&+ \underbrace{\psi_4\mathbb{G}\left(\pdiff{|A|^2}{\xi}\right) +\psi_5\mathbb{G}(\overline{A}B+A\overline{B})}_{\mathrm{Non-oscillating}} \qquad (\mbox{\textbf{Second-order solution}}),
\label{u3}
\end{split}
\eea
where $C(\xi,\tau)$ is an arbitrary modulation function and $\psi_i$ are functions of $k$.
\subsection{Two-time third-order: The modified-HNLS system}
The expansion at $O(\varepsilon^4)$ is
\begin{multline}
\mathcal{L}(u_4) + \pdiff{u_2}{\tau} + \frac{1}{2}\diff{c_g}{k}\pdiffn{u_2}{\xi}{2} + \frac{1}{6}\diffn{c_g}{k}{2}\pdiffn{u_1}{\xi}{3} + \mathcal{H}(u_1,u_3) + \mathcal{H}(u_3,u_1) +\mathcal{H}(u_2,u_2) \\+ \left[\diffslowoperator_{\mathcal{H},\xi}(u_1,u_2) + \diffslowoperator_{\mathcal{H},\xi}(u_2,u_1)\right] +\diffslowoperator_{\mathcal{H},\xi}^2(u_1,u_1)
+ \mathcal{T}(u_1,u_1,u_2) + \mathcal{T}(u_1,u_2,u_1)\\+\mathcal{T}(u_2,u_1,u_1) +\diffslowoperator_{\mathcal{T},\xi}(u_1,u_1,u_1) + \mathcal{Q}(u_1,u_1,u_1,u_1) + \mathcal{V}'\left(\pdiff{u_1}{\xi}\right)= 0,
\end{multline}
where $\mathcal{V}'$ is the Fre\'chet derivative of the operator $\mathcal{V}$. In Fourier representation this is
\begin{multline}
    \mathcal{L}(u_4) = -\underbrace{\left[\pdiff{B}{\tau}+ \frac{1}{2}\diff{c_g}{k}\pdiffn{B}{\xi}{2}\right.}_{\mbox{Order $\varepsilon$ NLS}} \underbrace{+ \mu_iB\mathbb{G}(|A|^2) + \mu_iA^2\overline{B} + \mu_i|A|^2B}_{\mbox{Coupling terms}}\\
   \underbrace{\left.+\frac{1}{6}\diffn{c_g}{k}{2}\pdiffn{A}{\xi}{3} + \mu_1|A|^2\pdiff{A}{\xi} + \mu_2 A^2\pdiff{\overline{A}}{\xi}+\mu_3 A\, \mathbb{G}\left(\pdiff{|A|^2}{\xi}\right) + \mathcal{V}'(k)\pdiff{A}{\xi}\right]\mbox{e}^{\mathrm{i}\Omega} }_{\mbox{Resonant}} \\  + \underbrace{\mu_4A^4\mbox{e}^{4\mathrm{i}\Omega} + \mu_5A^2\pdiff{A}{\xi}\mbox{e}^{3\mathrm{i}\Omega}}_{\mbox{Non-resonant}}+\underbrace{\left(\mu_6|A|^2A^2+\mu_{7}A^2\mathbb{G}(|A|^2)+\mu_{8}A\pdiffn{A}{\xi}{2}+\mu_{9}\left(\pdiff{A}{\xi}\right)^2\right)\mbox{e}^{2\mathrm{i}\Omega}}_{\mbox{Non-resonant}} \\\underbrace{+ \left(\mu_{10}|A|^4+ \mu_{11}|A|^2\mathbb{G}(|A|^2)+\mu_{12}\mathbb{G}(|A|^2)^2+\mu_{13}\pdiffn{|A|^2}{\xi}{2}\right)}_{\mbox{Non-resonant}} + \mbox{c.c},
    \label{thirdorder}
\end{multline}
where $\mu_i$ are functions of $k$ and the transport terms vanish due to the travelling wave frame of reference, \eqref{travellingwave}. We emphasise that the expressions on the first two lines in \eqref{thirdorder} are resonant.

Invoking the Fredholm alternative yields the solvability condition
\begin{multline}
\pdiff{B}{\tau}+ \frac{1}{2}\diff{c_g}{k}\pdiffn{B}{\xi}{2}+\mu_iB\mathbb{G}(|A|^2) + \mu_iA^2\overline{B} + +\mu_i|A|^2B \\
+\frac{1}{6}\diffn{c_g}{k}{2}\pdiffn{A}{\xi}{3} + \mu_1|A|^2\pdiff{A}{\xi} + \mu_2 A^2\pdiff{\overline{A}}{\xi} +\mu_3A\,\mathbb{G}\left(\pdiff{|A|^2}{\xi}\right) + \mathcal{V}'(k)\pdiff{A}{\xi} = 0.
\label{B_eqn}
\end{multline}
Equation \eqref{B_eqn} is an evolution equation for $B(\xi,\tau)$ that is coupled to $A(\xi,\tau)$ through \eqref{A_eqn} in the \textbf{modified-HNLS system}:
\begin{equation}
            \boxed{            \begin{aligned}
                                &\pdiff{A}{\tau} + \frac{1}{2}\diff{c_g}{k}\pdiffn{A}{\xi}{2} + \lambda_1 A|A|^2 + V(k)A+{\lambda_5A\mathbb{G}(|A|^2)} = 0, \\ &\pdiff{B}{\tau}+ \frac{1}{2}\diff{c_g}{k}\pdiffn{B}{\xi}{2}+\mu_iB\mathbb{G}(|A|^2) + \mu_iA^2\overline{B} + \mu_i|A|^2B \\
                                &+\frac{1}{6}\diffn{c_g}{k}{2}\pdiffn{A}{\xi}{3} + \mu_1|A|^2\pdiff{A}{\xi} + \mu_2 A^2\pdiff{\overline{A}}{\xi} +\mu_3A\,\mathbb{G}\left(\pdiff{|A|^2}{\xi}\right) + \mathcal{V}'(k)\pdiff{A}{\xi} = 0.
                \end{aligned}}
                \label{dysthe}
\end{equation}
After solving \eqref{dysthe} for $A$ and $B$ the asymptotically consistent solution to \eqref{pde1} is simply
\bea
u_{\varepsilon} = (A + \varepsilon B)\mbox{e}^{\mathrm{i}\Omega} + \varepsilon u_2,
\label{explicit}
\eea
where $u_2$ is stated in \eqref{u2}. Furthermore, we emphasise that \eqref{explicit} is the explicit second-order reconstruction of the solution $u$, once $A$ and $B$ have been found by solving \eqref{dysthe}, and is relatively straightforward to construct.

\subsection{Discussion}

We have derived a system of equations, \eqref{dysthe}, that determine an asymptotically consistent form of the wave-envelope. Crucially, the small parameter $\varepsilon$, is absent from the final form in \eqref{dysthe}, thus fulfilling our main aim. We emphasise that this analysis has been on a general PDE and has demonstrated the universality of the wave-envelope evolution equations in \eqref{dysthe} for systems defined in \eqref{pde1}. An important part of this procedure was Taylor expanding a nonlinear functional in terms of Fre\'chet derivatives. However, for specific systems knowledge of Fre\'chet derivatives is not needed as we shall now demonstrate for two specific systems; a toy fifth-order problem and the more complicated damped water-wave problem.

\section{Example 1: A toy fifth-order system}\label{sec:kdv}
We now apply this analysis to the toy dispersive system (without any damping)
\bea
\pdiff{u}{t}  + \pdiffn{u}{x}{3} + \pdiffn{u}{x}{5} + \varepsilon\left(u\pdiff{u}{x} + u\pdiffn{u}{x}{3} + \pdiff{u}{x}\pdiffn{u}{x}{2} \right)  = 0,\label{kdv}
\eea
which, with appropriately scaled coefficients is related to the fifth-order Korteweg-de Vries equation, see \cite{karima2018fifth}. We will expand \eqref{kdv} using \eqref{usol_best} with base-state $u_s = 0$. It is instructive to omit the multi-scale expansion at first to make the form of the operators transparent. Therefore, our expansion is 
\bea
u = \varepsilon u_1(x,t) + \varepsilon^2 u_2(x,t) + \varepsilon^3 u_3(x,t) + \varepsilon^4 u_4(x,t) + \cdots
\label{kdvexpansion}
\eea
Substituting \eqref{kdvexpansion} into \eqref{kdv} gives
\begin{multline}
  \varepsilon\mathcal{L}(u_1) + \varepsilon^2\left[\mathcal{L}(u_2) + \mathcal{H}(u_1,u_1)\right] + \varepsilon^3\left[\mathcal{L}(u_3) + \mathcal{H}(u_1,u_2) + \mathcal{H}(u_2,u_1) \right]\\
  + \varepsilon^4\left[\mathcal{L}(u_4) +\mathcal{H}(u_2,u_2) + \mathcal{H}(u_1,u_3) + \mathcal{H}(u_3,u_1) \right] + \cdots= 0,\label{kdvequationexpansion}
\end{multline}
where we have identified the operators
\begin{equation}
  \mathcal{L}(u) \equiv \pdiff{u}{t} + \mathcal{J}(u),\qquad \mathcal{J}(u)\equiv\pdiffn{u}{x}{3} + \pdiffn{u}{x}{5}, \qquad \mathcal{H}(u,v)\equiv u\pdiff{v}{x} + \pdiff{u}{x}\pdiffn{v}{x}{2} + u\pdiffn{v}{x}{3}.\label{kdvoperators}
\end{equation}
The Fourier multipliers of the operators in \eqref{kdvoperators} are
\bea
J(k_1) = \mathrm{i}k_1^3(k_1^2 - 1),\qquad H(k_1,k_2) = \mathrm{i}k_2(1 - k_1k_2 - k_2^2).
\label{kdvmultipliers}
\eea
We emphasise that we have \textit{not} found the Fre\'chet derivatives explicitly; the definitions in \eqref{kdvoperators} arise from simple algebra. We now expand the operators by simply plugging \eqref{kdvmultipliers} into the multi-scale expansions in \eqref{asd2}-\eqref{asd3}. Our problem then becomes
\begin{multline}
  \varepsilon\mathcal{L}(u_1) + \varepsilon^2\left[\mathcal{L}(u_2) + \mathcal{H}(u_1,u_1) + \diffslowoperator_{\mathcal{J}}(u_1) + \mathrm{i}\pdiff{u_1}{T}\right] \\
  + \varepsilon^3\left[\mathcal{L}(u_3) + \mathcal{H}(u_1,u_2) + \mathcal{H}(u_2,u_1) + \diffslowoperator_{\mathcal{J}}^2(u_1) + \diffslowoperator_{\mathcal{H}}(u_1,u_1) + \diffslowoperator_{\mathcal{J}}(u_2) + \mathrm{i}\pdiff{u_2}{T}\right]\\
  + \varepsilon^4\left[\mathcal{L}(u_4) +\mathcal{H}(u_2,u_2) + \mathcal{H}(u_1,u_3) + \mathcal{H}(u_3,u_1) + \diffslowoperator_{\mathcal{J}}^3(u_1) + \diffslowoperator_{\mathcal{H}}^2(u_1,u_1) + \right.\\\left.
    \diffslowoperator_{\mathcal{J}}^2(u_2) + \diffslowoperator_{\mathcal{H}}(u_2,u_2) +  \diffslowoperator_{\mathcal{H}}(u_1,u_3) +  \diffslowoperator_{\mathcal{H}}(u_3,u_1) + \mathrm{i}\pdiff{u_3}{T}\right] + \cdots= 0.\label{kdvequationexpansion1}
\end{multline}
We are now in a position to simply `pick' out the equations at each order.
\subsection{First-order}
The expansion at $O(\varepsilon)$ \eqref{kdvequationexpansion1} is
\bea
\mathcal{L}(u_1) = 0,
\eea
which has the solution and dispersion relation
\bea
u_1 = A(X,T)\,\mbox{e}^{\mathrm{i}\Omega} + \mbox{c.c}, \quad \omega - k^3 + k^5= 0.
\label{kdvsol1}
\eea
\subsection{Second-order}
Continuing, at $O(\varepsilon^2)$ \eqref{kdvequationexpansion1} is
\bea
\mathcal{L}(u_2) = \mathrm{i}\left[\pdiff{A}{T}  + c_g\pdiff{A}{X}\right]\,\mbox{e}^{\mathrm{i}\Omega} + \mathrm{i}k(1 - 2k^2) A^2 \mbox{e}^{2\mathrm{i}\Omega},\qquad c_g = 5k^4 - 3k^2.
\label{kdvordereps}
\eea
Notice there are no non-oscillating terms in \eqref{kdvordereps}. We eliminate the resonant transport terms by moving to a frame of reference moving with the group velocity $c_g$ by setting $\xi = X - c_gT$ and $\tau = \varepsilon^2 t$. By finding the particular integral of \eqref{kdvordereps} and using the dispersion relation in \eqref{kdvsol1}, the solution at this order is therefore
\bea
u_2 = B(\xi,\tau)\,\mbox{e}^{\mathrm{i}\Omega} + \varphi_0A^2\mbox{e}^{2\mathrm{i}\Omega} + \mbox{c.c.}, \quad \varphi_0 = \frac{1-2k^2}{2k^2(17k-5)}.
\eea
\subsection{Third-order} 
The next order in \eqref{kdvequationexpansion1} is $O(\varepsilon^3)$ where 
\begin{multline}
\mathcal{L}(u_3) = -\mathrm{i}\left[\mathrm{i}\pdiff{A}{\tau}  + k(1-5k^4)A|A|^2 + \diff{c_g}{k}\pdiffn{A}{\xi}{2}\right]\,\mbox{e}^{\mathrm{i}\Omega} \\+ \left((2k^2 - 3k - 1)\pdiff{A}{\xi}A+k(2k^2-1)AB\right)\mbox{e}^{2\mathrm{i}\Omega}- \mathrm{i}\,4k(1 - 6k^2)\varphi_0A^3\mbox{e}^{3\mathrm{i}\Omega}.
\label{kdvordereps2}
\end{multline}
Using the Fredholm alternative we get the modified-NLS equation for our toy system
\bea
\boxed{\mathrm{i}\pdiff{A}{\tau}  + k(1-5k^4)A|A|^2 + k(20k^2 - 9)\pdiffn{A}{\xi}{2} = 0.}
\label{nlskdv}
\eea
The solution at this order is
\bea
u_3 = C(\xi,\tau)\,\mbox{e}^{\mathrm{i}\Omega} + \psi_1A\pdiff{A}{\xi}\mbox{e}^{2\mathrm{i}\Omega} +\psi_2AB\mbox{e}^{2\mathrm{i}\Omega} + \psi_3\mbox{e}^{3\mathrm{i}\Omega} + \mbox{c.c},
\eea
where $C(\xi,\tau)$ is undetermined at this order and by finding the particular integral of \eqref{kdvordereps2} and using the dispersion relation in \eqref{kdvsol1}:
\bea
\psi_1 = \frac{2k^2-1}{2k(17k-5)},\qquad \psi_2 = \frac{2k^2 - 3k - 1}{2k^3(17k-5)},\qquad \psi_3 = \frac{2(1-6k^2)\varphi_0}{3k^2(41k^2 - 5)}.
\eea
\subsection{Fourth-order}
Finally, at $O(\varepsilon^4)$ \eqref{kdvequationexpansion1} is
\begin{multline}
\mathcal{L}(u_4) = \mathrm{i}\left[\mathrm{i}\pdiff{B}{\tau} - k(1-5k^4)B|B|^2 + \diff{c_g}{k}\pdiffn{B}{\xi}{2} + 2\psi_1|A|^2\pdiff{A}{\xi} - 2\psi_2|A|^2 B - \varphi_0A^2\overline{B} \right.\\\left.+ (1+4k)A^2\pdiff{\overline{A}}{\xi} + (1 - 4k + 3k^2)\overline{A}\pdiff{|A|^2}{\xi} + 6k(10k-3)\pdiffn{A}{X}{3}\right]\mbox{e}^{\mathrm{i}\Omega} + \mbox{Non-resonant terms}.
\label{kdvordereps3}
\end{multline}
Using the Fredholm alternative we get the evolution equation for $B$:
\begin{multline}
  \mathrm{i}\pdiff{B}{\tau} - k(1-5k^4)B|B|^2 + \diff{c_g}{k}\pdiffn{B}{\xi}{2} + 2\psi_1|A|^2\pdiff{A}{\xi} - 2\psi_2|A|^2 B - \varphi_0A^2\overline{B} \\+ (1+4k)A^2\pdiff{\overline{A}}{\xi} + (1 - 4k + 3k^2)\overline{A}\pdiff{|A|^2}{\xi} + 6k(10k-3)\pdiffn{A}{X}{3} = 0\label{hnlskdv}
\end{multline}
Equation \eqref{hnlskdv} coupled with \eqref{nlskdv} form the complete HNLS for the toy problem in \eqref{kdv}:
\bea
\boxed{\begin{aligned}
    &\mathrm{i}\pdiff{A}{\tau}  + k(1-5k^4)A|A|^2 + k(20k^2 - 9)\pdiffn{A}{\xi}{2} = 0\vspace{1cm}\\
    &\mathrm{i}\pdiff{B}{\tau} - k(1-5k^4)B|B|^2 + \diff{c_g}{k}\pdiffn{B}{\xi}{2} + 2\psi_1|A|^2\pdiff{A}{\xi} - 2\psi_2|A|^2 B - \varphi_0A^2\overline{B} \\
    &+ (1+4k)A^2\pdiff{\overline{A}}{\xi} + (1 - 4k + 3k^2)\overline{A}\pdiff{|A|^2}{\xi} + 6k(10k-3)\pdiffn{A}{X}{3} = 0
		\label{dysthekdv}
\end{aligned}}
\eea

\subsection{A note on the boundary conditions}

As a final note, \eqref{kdv} is fifth-order in space and so five conditions on the $x-$derivatives of $u$ are required to ensure the problem is well-posed. The system in \eqref{dysthekdv} is third-order in space and hence requires three boundary conditions and the other two are accounted for by the periodicity requirement. Therefore, going to higher-orders for this toy system is ill-posed as more boundary conditions for the wave envelope are required than the original toy problem in \eqref{dysthekdv}

\section{Example 2: The damped water-wave problem}\label{sec:water}

We now describe a more complicated example; determining the shape of the free-surface, $y = \zeta$, bounding a body of irrotational, inviscid and incompressible fluid with damping. In the fluid, the unknown velocity potential, $\phi$, must satisfy Laplace's equation in the fluid, and in this instance, we neglect surface tension. 

 We nondimensionalise the system in the same way as \cite{lannes2013book}. To mimic the general system in \eqref{pde1}, the governing equations can be written as a dynamical system of the form
\bea
\pdiff{\boldsymbol{u}}{t} + \boldsymbol{\mathcal{F}}(\boldsymbol{u}) + \varepsilon^2\boldsymbol{\mathcal{V}}(\boldsymbol{u})=0,
\eea
where $\boldsymbol{u} = [\zeta,\psi ]^T$, $\psi  = \phi_{y=\zeta}$ and $\boldsymbol{\mathcal{V}}$ is the, now vectorial, linear dissipation term which is consistent with the approach of \cite{dias2008pla}. The nonlinear functional can be written explicitly as 
\bea
\boldsymbol{\mathcal{F}}(\boldsymbol{u}) = \left[-\mathcal{G}[\varepsilon\zeta]\psi ,\zeta +\frac{1}{2}\varepsilon|\nabla\psi |^2 -\frac{\varepsilon(\mathcal{G}[\varepsilon\zeta](\psi ) + \varepsilon\nabla\zeta\cdot\nabla\psi )^2}{2(1+\varepsilon^2|\nabla\zeta|^2)}\right]^T,\label{water_wave_dyn}
\eea
where $\varepsilon$ is the nonlinearity parameter, proportional to the ratio of a typical wave-height to wave-length. The spatial boundary conditions are that $\varphi$ and $\zeta$ and the first spatial derivative of $\zeta$ vanish as $x\to\pm\infty$. In addition, $\nabla = [\partial/\partial x,\partial/\partial y]^T$, $G[\zeta]$ is the Dirichlet-to-Neumann map (DtN) on deep water defined as
\bea
\mathcal{G}[\zeta](\psi ) = (1 + |\nabla\zeta|^2)^{1/2}\left(\nabla\phi\cdot\boldsymbol{n}\right),
\label{dtn_first}
\eea
where $\boldsymbol{n}$ is the outwards-pointing normal of the free-surface. The operator is linear with respect to $\psi $ (and hence the Fre\'chet derivatives with respect to $\psi $ are trivial) but highly nonlinear with respect to $\zeta$ and more work is required to determine the Fre\'chet derivatives with respect to $\zeta$.

As the dynamical system is two dimensional in the sense there are two independent phase variables in $\boldsymbol{u}$, the coefficients of the modified-HNLS equation for this system are vectors in $\mathbb{R}^2$ and the envelope is also a vector, $\boldsymbol{A} =[A_{\zeta},A_{\psi }]^T$. Furthermore, the domain is spatially two-dimensional so, to reflect \eqref{usol_best}, we will look for a solution of the form
\bea
\boldsymbol{u} = \boldsymbol{A}\,\mbox{e}^{\mathrm{i}(\boldsymbol{k}\cdot\boldsymbol{x}-\omega t)}+ \varepsilon\boldsymbol{u}_2 + \varepsilon^2\boldsymbol{u}_3 + \cdots + \mbox{c.c},
\eea
where $\boldsymbol{u}_n = [\zeta_n,\psi _n]^T$ and $\boldsymbol{k} = [k_x,k_y]^T$.

A key part of the analysis is the expansion of the DtN in \eqref{dtn_first}. The DtN is analytic in $\zeta$ (see, for example \cite{hu2005analyticity}) and hence we expand it as a Taylor series about the zero base state
 \bea
  \mathcal{G}[\varepsilon\zeta](\psi ) = \sum_{n=0}^{\infty}\varepsilon^n\mathcal{G}_n[\zeta](\psi )
\label{dtn}
 \eea
where $\mathcal{G}_n$ is an $n$-linear operator. The form of $\mathcal{G}_n$ is highly dependent on the \textit{geometry} of the problem. For a finite depth domain with a flat bottom and for an infinite depth domain the explicit forms of $\mathcal{G}_n$ are known, for example \cite{wilkening2015comparison} and \cite{guyenne2022ham}, respectively. Whilst the DtN is well defined for more complex geometries, including varying bottom topography, the explicit calculation is more difficult \cite{craig2005hamiltonian,andrade2018three}. As the main focus of this article is to show that the evolution equations emerge from general dispersive systems with weak dissipation, not just for water-wave problems, we do not explore these difficulties which are well documented in the literature. With that in mind, we apply our methodology to waves over infinite depth, in which case the expansion in \eqref{dtn} is relatively straightforward. For infinite-depth, the first few terms of the DtN are
\begin{align}
&\mathcal{G}_0(\psi) = |\diffoperator|(\psi),\qquad \mathcal{G}_1[\zeta](\psi) = -|\diffoperator|(\zeta(|\diffoperator|(\psi)) - \nabla\cdot(\zeta\nabla\psi),\\ &\mathcal{G}_2[\zeta](\psi) = |\diffoperator|(\zeta|\diffoperator|(\zeta|\diffoperator|(\psi))) + \frac{1}{2}\nabla^2(\zeta^2|\diffoperator|(\psi)) + \frac{1}{2}|\diffoperator|(\zeta^2\nabla^2\psi),
\end{align}
where $\diffoperator=-\mathrm{i}\,\nabla$ with a Fourier multiplier of $\boldsymbol{k}$, and $|\diffoperator|$ has a Fourier multiplier of $|\boldsymbol{k}|$, see p.88 of \cite{lannes2013book}. 

\subsection{Leading-order}
At $O(1)$ we obtain 
\bea
\boldsymbol{\mathcal{L}}(\boldsymbol{u}_1)\equiv \pdiff{\boldsymbol{u}_1}{t} + \boldsymbol{\mathcal{J}}(\boldsymbol{u}_1)=0,\qquad \boldsymbol{\mathcal{J}}= \begin{pmatrix}
0 & -\mathcal{G}_0\\
1 & 0
\end{pmatrix},\qquad \boldsymbol{J}(\boldsymbol{k}_1)=\begin{pmatrix}
0 & -|\boldsymbol{k}|\\
1 & 0
\end{pmatrix},
\eea
where the bold symbol are matrix versions of $\mathcal{J} (u_1)$ and $J(k)$ introduced earlier. Writing $\boldsymbol{u}_1$ in component form,
\bea
\boldsymbol{u}_1 =[\zeta_1, \psi _1]^T= \boldsymbol{A}\mbox{e}^{\mathrm{i}(\boldsymbol{k}\cdot\boldsymbol{x}-\omega t)},
\eea
results in an eigenvalue problem. To ensure that a non-zero solution exists the full dispersion relation for deep-water waves emerges
\bea
\omega =J(\boldsymbol{k})\equiv |\boldsymbol{k}|^{1/2}.
\label{dispersionwater}
\eea
and 
\bea
\boldsymbol{A}=A(\boldsymbol{X},T)
\boldsymbol{v},\qquad \boldsymbol{v} = [\mathrm{i}\omega,1]^T,
\eea
where $\boldsymbol{X} = \varepsilon\boldsymbol{x}$ and $T=\varepsilon t$. Note the dispersion relation \eqref{dispersionwater} does not contain the acceleration due to gravity (as one might expect) due to the nondimensionalisation. 

It is also useful to Taylor expand the dispersion relation and the operator $\mathcal{G}_0$ in the slow-variable. We find that
\begin{align}
	&\omega(\boldsymbol{k})\mapsto \omega(\boldsymbol{k}) - \varepsilon\mathrm{i}\nabla_{\boldsymbol{k}}\omega\cdot\nabla_{\boldsymbol{X}} - \frac{1}{8}\varepsilon^2\frac{\nabla_{\boldsymbol{X}}^2}{|\boldsymbol{k}|^{3/2}} + \cdots \\
	&\mathcal{G}_0(A(\boldsymbol{\xi})\mbox{e}^{\mathrm{i}\Omega})\mapsto A(\boldsymbol{\xi})|\boldsymbol{k}|\mbox{e}^{\mathrm{i}\Omega} - \varepsilon\left(2\mathrm{i}\omega\nabla_{\boldsymbol{k}}\omega(\boldsymbol{k})\cdot \nabla_{\boldsymbol{\xi}}A\right) + O(\varepsilon^3),
\end{align}
where $\nabla_{\boldsymbol{k}} = [\partial/\partial k_1,\partial/\partial k_2]^T$ and $\nabla_{\boldsymbol{\xi}} = [\partial/\partial \xi_1,\partial/\partial \xi_2]^T$. We also note that
\bea
\nabla_{\boldsymbol{k}}\omega = \frac{\boldsymbol{k}}{2\omega^3},\qquad \nabla^2_{\boldsymbol{k}}\omega = -\frac{1}{4\omega^3}
\eea
Finally, a particular quirk of \eqref{water_wave_dyn} is that $\phi$ only appears as a derivative, and hence we are free to add an arbitrary correction to the velocity potential which, as we shall see later, is crucial for the determination of the evolution equations. Therefore, the first-order solution is 
\bea
\boldsymbol{u}_1 = A\boldsymbol{v}\mbox{e}^{\mathrm{i}\Omega} + [0,\bar{\psi}_1]^T \qquad (\mbox{\textbf{First-order solution}}),
\eea
where $\bar{\psi}_1$ is arbitrary at this stage and represents a mean-flow.

\subsection{First-order}
The expansion at $O(\varepsilon)$ is
\bea
\boldsymbol{\mathcal{L}}(\boldsymbol{u}_2)= 
\underbrace{\left[\pdiff{A}{T} + 
\nabla\omega\cdot(\nabla_{\boldsymbol{X}}A)
\right]}_{\mbox{Resonant terms}}\boldsymbol{v}\,\mbox{e}^{\mathrm{i}\Omega} + \mbox{c.c}
    +
    \boldsymbol{\mathcal{H}}(\boldsymbol{u}_1,\boldsymbol{u}_1),
 \eea
 where
 \bea
 \boldsymbol{\mathcal{H}}(\boldsymbol{u}_i,\boldsymbol{u}_j) \equiv 
 \begin{pmatrix}
        \mathcal{G}_1[\zeta_i](\psi_j)\\
	 -\frac{1}{2}\nabla\psi_i\cdot\nabla\psi_j +\frac{1}{2}\mathcal{G}_0(\psi_i)\mathcal{G}_0(\psi_j)
     \end{pmatrix}.
\eea
As before, to eliminate the secular terms above we set
\bea
\boldsymbol{\xi} = \boldsymbol{X} - \boldsymbol{c}_g T\qquad\boldsymbol{c}_g = \nabla_{\boldsymbol{k}}\omega.\label{water_wave_frame}
\eea

The solution to $\boldsymbol{u}_2$ is
\bea
\boldsymbol{u}_2 = \boldsymbol{\varphi}_0\,A^2\mbox{e}^{2\mathrm{i}\Omega}  + B\,\boldsymbol{v}\mbox{e}^{\mathrm{i}\Omega} + \mbox{c.c} + [0,\bar{\varphi}_2]^T,
\eea
where $B(\boldsymbol{\xi},T)$ and $\bar{\varphi}_2(\boldsymbol{\xi})$ are both arbitrary at this order and $\boldsymbol{\varphi}_0$ solves the linear problem
\bea
\boldsymbol{L}_2\boldsymbol{\varphi}_0 = [    0,\,
    |\boldsymbol{k}|^2]^T, \qquad \boldsymbol{L}_n \equiv \begin{pmatrix}
    -n\mathrm{i}\omega & -n|\boldsymbol{k}|\\
    1 & -n\mathrm{i}\omega
\end{pmatrix}.
\label{varphi_solution}
\eea
yielding
\bea
\boldsymbol{\varphi}_0 = \frac{\mathrm{i}|\boldsymbol{k}|^2}{\omega}\,\boldsymbol{v}.
           \eea
We conclude that at $O(\varepsilon)$ there are no non-oscillating terms in $\boldsymbol{u}_2$. 

\subsection{Second-order}
Continuing the expansion, at $O(\varepsilon^2)$ (by moving to a frame of reference using \eqref{water_wave_frame} and using the small timescale $\tau = \varepsilon^2t$) we obtain
\begin{multline}
\boldsymbol{\mathcal{L}}(\boldsymbol{u}_3)= 
-\left[\pdiff{A}{\tau} - \frac{1}{2}\mathrm{i}\nabla^2\omega
	\nabla^2_{\boldsymbol{\xi}}A\right]\boldsymbol{v}\,\mbox{e}^{\mathrm{i}\Omega} + \mbox{c.c} + \diffslowoperator_{\boldsymbol{\mathcal{H}}}(\boldsymbol{u}_1,\boldsymbol{u}_1) + \boldsymbol{\mathcal{H}}(\boldsymbol{u}_1,\boldsymbol{u}_2) + \boldsymbol{\mathcal{H}}(\boldsymbol{u}_2,\boldsymbol{u}_1)\\
	+\boldsymbol{\mathcal{T}}(\boldsymbol{u}_1,\boldsymbol{u}_1,\boldsymbol{u}_1),
\label{bethal}
\end{multline}
where
\bea
\boldsymbol{\mathcal{T}}(\boldsymbol{u}_i,\boldsymbol{u}_j,\boldsymbol{u}_k) \equiv
	\begin{pmatrix}
    \mathcal{G}_2[\zeta_i,\zeta_j](\psi_k)\\
\mathcal{G}_0(\psi_i)\mathcal{G}_1[\zeta_j](\psi_k) 
+\mathcal{G}_0(\psi_i)\nabla\zeta_j\cdot\nabla\psi_k.
\end{pmatrix}
\eea
In Fourier space \eqref{bethal} is
\begin{multline}
\boldsymbol{\mathcal{L}}(\boldsymbol{u}_3)= 
	-\underbrace{\left[\pdiff{A}{\tau} - \frac{1}{2}\mathrm{i}\nabla^2\omega\nabla^2_{\boldsymbol{\xi}}A \right]\boldsymbol{v}\,\mbox{e}^{\mathrm{i}\Omega} - \frac{\mathrm{i}|\boldsymbol{k}|^4}{\omega} A|A|^2\begin{pmatrix}3\mathrm{i}\omega\\ 1\end{pmatrix}\mbox{e}^{\mathrm{i}\Omega}}_{\mbox{Resonant terms}} \\\underbrace{+ \mathcal{\boldsymbol{{V}(\boldsymbol{k})}}A\,\boldsymbol{v}\mbox{e}^{\mathrm{i}\Omega} + \diffslowoperator_{\boldsymbol{\mathcal{H}}}([0,\bar{\psi}_1]^T,\boldsymbol{u}_1) +\diffslowoperator_{\boldsymbol{\mathcal{H}}}(\boldsymbol{u}_1,[0,\bar{\psi}_1]^T)}_{\mbox{Resonant terms}} \\+\underbrace{\boldsymbol{\lambda}_1 AB
\mbox{e}^{2\mathrm{i}\Omega} + \boldsymbol{\lambda}_2 A^3\mbox{e}^{3\mathrm{i}\Omega} + \mbox{c.c} + \boldsymbol{\lambda}_3(A\overline{B}
- \overline{A}B)
}_{\mbox{Non-resonant terms}},
	\label{fourier_space}
\end{multline}
where
\bea
\boldsymbol{\lambda}_1 = [0,\,2|\boldsymbol{k}|^2]^T,\qquad \boldsymbol{\lambda}_2 = [0,\,3\omega|\boldsymbol{k}|^3]^T,\qquad \boldsymbol{\lambda}_3 = [0,\,|\boldsymbol{k}|^2]^T,
\eea
and $\boldsymbol{V}(\boldsymbol{k})$ is the pseudo-operator of $\boldsymbol{\mathcal{V}}$. At this stage, the resonant terms are not consistent as they do not contain a common factor of $\boldsymbol{v}$. This can be rectified by a judicious choice of our arbitrary leading-order mean flow, $\bar{\psi}_1$. The fourth and fifth resonant term in \eqref{fourier_space} can be written explicitly as 
\bea
\diffslowoperator_{\boldsymbol{\mathcal{H}}}([0,\bar{\psi}_1]^T,\boldsymbol{u}_1) +\diffslowoperator_{\boldsymbol{\mathcal{H}}}([\bar{\psi}_1,0]^T,\boldsymbol{u}_1) = \begin{pmatrix}
	-|\diffoperator|\left(\zeta_1 \frac{\diffoperator}{|\diffoperator|}\,\diffoperator_{\boldsymbol{\xi}}(\bar{\psi}_1)\right) + \diffoperator(\zeta_1\diffoperator_{\boldsymbol{\xi}}(\bar{\psi}_1))\\
	\diffoperator_{\boldsymbol{\xi}}\bar{\psi}_1\cdot\nabla\psi_1 + \frac{|\diffoperator|}{\diffoperator}\,\diffoperator_{\boldsymbol{\xi}}\bar{\psi}_1|\diffoperator|(\psi_1)
\end{pmatrix}.
\label{hessian_mean_flow}
\eea
By choosing the $\psi$ component in \eqref{hessian_mean_flow} to vanish, we obtain the desired action of the non-local functional, $\diffoperator/|\diffoperator|$, acting on $\diffoperator_{\mathcal{H}}\bar{\psi}_1$, i.e.
\bea
\frac{\diffoperator}{|\diffoperator|}\diffoperator_{\boldsymbol{\xi}}\bar{\psi}_1 = -\frac{\boldsymbol{k}}{|\boldsymbol{k}|}\diffoperator_{\boldsymbol{\xi}}\bar{\psi}_1.
\label{useful_relation}
\eea
Substituting \eqref{useful_relation} into the $\zeta$ component of \eqref{hessian_mean_flow} and insisting that the contribution adds to the second resonant term in \eqref{fourier_space}, in such a way that a factor of $\boldsymbol{v}$ can be factored out, means that we require $\bar{\psi}_1$ to satisfy
\bea
\boldsymbol{k}\cdot\nabla_{\boldsymbol{\xi}}\bar{\psi}_1 = -\frac{|\boldsymbol{k}|^4}{\omega}|A|^2.
\label{mean_flow}
\eea
With this choice of $\bar{\psi}_1$ the resonant terms in \eqref{fourier_space} all contain a common factor of $\boldsymbol{v}$ and hence by imposing the Fredholm alternative we arrive at the modified-NLS equation:
\bea
\boxed{\mathrm{i}\pdiff{A}{\tau} + \left(\frac{1}{2}\nabla^2\omega\right)\nabla^2_{\boldsymbol{\xi}}A -\left(\frac{|\boldsymbol{k}|^4}{\omega}\right)A|A|^2 + \boldsymbol{V}(\boldsymbol{k})A = 0.}
\label{nlswaves}
\eea
Notice that for this particular example of a dynamical system, no non-local terms are entering the modulation equation, unlike the general modified-NLS in \eqref{A_eqn}.

Continuing further, the solution at $O(\varepsilon^2)$ is
\bea
    \boldsymbol{u}_3 = \underbrace{\boldsymbol{\varphi}_2A^3\boldsymbol{v}\mbox{e}^{3\mathrm{i}\Omega}+ \left(\boldsymbol{\varphi}_3AB +(\boldsymbol{\varphi}_{4}\cdot A\nabla_{\boldsymbol{\xi}}A)\right)\boldsymbol{v}\mbox{e}^{2\mathrm{i}\Omega}+ C\boldsymbol{v}\mbox{e}^{\mathrm{i}\Omega} + \mbox{c.c}}_{\mbox{Oscillating terms}} + \underbrace{\begin{pmatrix}
        0\\
        \frac{2\omega\boldsymbol{k}\cdot\nabla_{\boldsymbol{\xi}}|A|^2}{|D|}
    \end{pmatrix} }_{\mbox{Non-oscillating terms}},
\eea
where
\bea
\boldsymbol{\varphi}_2 = 3|\boldsymbol{k}|^3,\qquad \boldsymbol{\varphi}_3 = \frac{2\mathrm{i}|\boldsymbol{k}|^2}{\omega}\boldsymbol{v},\qquad \boldsymbol{\varphi}_4 = \frac{2\boldsymbol{k}}{\omega}.
\eea

\subsection{Third-order}
At $O(\varepsilon^3)$ the linear system starts has a large number of terms but we just concentrate on those terms that are resonant and hence contribute to the imposition of the Fredholm alternative:
\begin{multline}
	\boldsymbol{\mathcal{L}}(\boldsymbol{u}_4) = -\left[\pdiff{B}{\tau} -\frac{1}{2}\mathrm{i}\nabla^2\omega\nabla_{\boldsymbol{\xi}}^2 B\right]\,\boldsymbol{v}\mbox{e}^{\mathrm{i}\Omega} + \diffslowoperator_{\boldsymbol{\mathcal{H}}}(\boldsymbol{u}_1,\boldsymbol{u}_2)+ \diffslowoperator_{\boldsymbol{\mathcal{H}}}(\boldsymbol{u}_2,\boldsymbol{u}_1) + \diffslowoperator_{\boldsymbol{\mathcal{T}}}(\boldsymbol{u}_1,\boldsymbol{u}_1,\boldsymbol{u}_1) +\\
	\boldsymbol{\mathcal{H}}(\boldsymbol{u}_1,\boldsymbol{u}_3) +\boldsymbol{\mathcal{H}}(\boldsymbol{u}_3,\boldsymbol{u}_1) +  \boldsymbol{\mathcal{H}}(\boldsymbol{u}_2,\boldsymbol{u}_2) +  	\boldsymbol{\mathcal{T}}(\boldsymbol{u}_1,\boldsymbol{u}_1,\boldsymbol{u}_2) +\boldsymbol{\mathcal{T}}(\boldsymbol{u}_1,\boldsymbol{u}_2,\boldsymbol{u}_1) +\boldsymbol{\mathcal{T}}(\boldsymbol{u}_2,\boldsymbol{u}_1,\boldsymbol{u}_1)\\
	+ \mbox{Non-resonant terms}.	
\end{multline}
As in the previous order, to ensure the resonant terms contain a common factor of $\boldsymbol{v}$ requires that
\bea
\frac{\diffoperator}{|\diffoperator|}\diffoperator_{\boldsymbol{\xi}}(\bar{\psi}_2) = -\frac{\boldsymbol{k}}{|\boldsymbol{k}|}\diffoperator_{\boldsymbol{\xi}}(\bar{\psi}_2),
\eea
and 
\bea
\boldsymbol{k}\cdot \diffoperator_{\boldsymbol{\xi}}(\bar{\psi}_2) = \frac{\mathrm{i}|\boldsymbol{k}|^4}{\omega}\overline{A}B + \frac{\mathrm{i}|\boldsymbol{k}|^4}{2\omega}A\overline{B} - \frac{5|\boldsymbol{k}|^2}{\omega}\boldsymbol{k}\cdot \overline{A}\nabla_{\boldsymbol{\xi}}A - \frac{|\boldsymbol{k}|^2}{\omega}\boldsymbol{k}\cdot A\nabla_{\boldsymbol{\xi}}\overline{A} + \frac{2B}{A}|\boldsymbol{k}|\frac{\diffoperator}{|\diffoperator|}\diffoperator_{\boldsymbol{\xi}}(\bar{\psi}_1).
\eea
Hence the $O(\varepsilon^3)$ solvability condition yields
\bea
{\begin{aligned}
	\mathrm{i}\pdiff{B}{\tau} - \frac{1}{2}\mathrm{i}\nabla^2\omega\nabla_{\boldsymbol{\xi}}^2B - \frac{\mathrm{i}|\boldsymbol{k}|^2}{\omega}|A|^2B - \frac{\mathrm{i}|\boldsymbol{k}|^2}{\omega}A^2\overline{B} + \frac{2|\boldsymbol{k}|^2}{\omega}\boldsymbol{k}\cdot |A|^2\nabla_{\boldsymbol{\xi}}A + \frac{|\boldsymbol{k}|^2}{\omega}\boldsymbol{k}\cdot A^2\nabla_{\boldsymbol{\xi}}\overline{A}\\
	+ \frac{1}{24}\frac{\boldsymbol{k}}{\omega^7}\cdot(\nabla_{\boldsymbol{\xi}}\nabla^2_{\boldsymbol{\xi}}A) +	\underbrace{\frac{|\boldsymbol{k}|^3}{\omega}\frac{\diffoperator}{|\diffoperator|}\nabla_{\boldsymbol{\xi}}|A|^2 + B|\boldsymbol{k}|\frac{\diffoperator}{|\diffoperator|}(\diffoperator_{\boldsymbol{\xi}}(\bar{\psi}_1))}_{\mbox{Non-local terms}} = 0.
		\label{dysthewaves}
\end{aligned}}
\eea
Finally, we write down our coupled higher-order system for the slow modulation amplitudes, $A$ and $B$:
\bea
\boxed{\begin{aligned}
	&\mathrm{i}\pdiff{A}{\tau} + \left(\frac{1}{2}\nabla^2\omega\right)\nabla^2_{\boldsymbol{\xi}}A -\left(\frac{|\boldsymbol{k}|^4}{\omega}\right)A|A|^2 + \boldsymbol{V}(\boldsymbol{k})A=0. \\\vspace{1cm}
		&\mathrm{i}\pdiff{B}{\tau} - \frac{1}{2}\mathrm{i}\nabla^2\omega\nabla_{\boldsymbol{\xi}}^2B - \frac{\mathrm{i}|\boldsymbol{k}|^2}{\omega}|A|^2B - \frac{\mathrm{i}|\boldsymbol{k}|^2}{\omega}A^2\overline{B} + \frac{2|\boldsymbol{k}|^2}{\omega}\boldsymbol{k}\cdot |A|^2\nabla_{\boldsymbol{\xi}}A + \frac{|\boldsymbol{k}|^2}{\omega}\boldsymbol{k}\cdot A^2\nabla_{\boldsymbol{\xi}}\overline{A}\\
    &+	{ \frac{1}{24}\frac{\boldsymbol{k}}{\omega^7}\cdot(\nabla_{\boldsymbol{\xi}}\nabla^2_{\boldsymbol{\xi}}A) + \frac{|\boldsymbol{k}|^3}{\omega}\frac{\diffoperator}{|\diffoperator|}\nabla_{\boldsymbol{\xi}}|A|^2 + B|\boldsymbol{k}|\frac{\diffoperator}{|\diffoperator|}(\diffoperator_{\boldsymbol{\xi}}(\bar{\psi}_1))} + \nabla_{\boldsymbol{\xi}}A\cdot \diffoperator_{\boldsymbol{\xi}}(\boldsymbol{{V}})= 0,\\
    &\hspace{4.5cm}\boldsymbol{k}\cdot\nabla_{\boldsymbol{\xi}}\bar{\psi}_1 = -\frac{|\boldsymbol{k}|^4}{\omega}|A|^2.
		\label{dysthewaves1}
\end{aligned}}
\eea
We note a feature of this system is that the mean flow term, $\bar{\psi}_1$ is coupled to the first-order wave envelope via \eqref{mean_flow} a feature shared with the Davey-Stewartson system; see, \cite{davey1974,lannes2013book}.

\section{Discussion}\label{sec:discussion}
In this article, we have formally derived the modified-NLS and -HNLS systems for a general set of dispersive PDEs with damping. We have achieved this by posing the system as an infinite-dimensional dynamical system and performing a series of Taylor expansions on the nonlinear functionals, allowing us to derive general evolution equations using the method of multiple scales. In particular, we have placed the recent work of \cite{alberello2022dissipative,alberello2023dynamics} on a firmer theoretical footing. In \cite{alberello2022dissipative,alberello2023dynamics} a modified-NLS equation damping was used to describe ocean waves in the marginal ice-zone of the southern ocean.

We emphasise that the derivation of the evolution equations for this problem has been achieved before by many different authors (see introduction). Yet, we believe this general analytic framework provides a fresh perspective in that; i) it makes no assumption on what the small parameter, $\varepsilon$ is; ii) there is no assumption that the underlying base-state is the zero state; iii) there is no requirement that the particular form of solutions should be known \textit{a priori}, as is often the case in existing derivations; and iv) this framework can easily be applied to similar problems, including those where a free-surface is absent in the formulation.

Although the evolution equations have been traditionally derived in the context of measuring the wave envelope of water-waves and optics, there is an intriguing possibility that the framework we have derived here can be used to approximate general time-dependent periodic invariant solutions of infinite-dimensional PDE in numerical formulations. An attractive feature of this formulation is that in any numerical computation of these systems, we can approximate the multi-linear operators using the numerically-constructed Jacobian and Hessian operators thus avoiding the complicated algebra involved in more complicated systems. This substantially widens the potential scope of our approach to a wider class of problems. For example, a recently discovered periodic invariant-solution of the water-wave problem for flow over a localised topographic forcing has been discovered numerically \cite{keeler2024hydraulic} and it is an intriguing possibility that the envelope equations developed here can be used to approximate such objects.

\vskip2pc


\bibliographystyle{jfm}
\bibliography{nls.bib}

\end{document}